\documentclass[
    prd,
    amssymb,
    preprintnumbers,
    secnumarabic,
    nofootinbib,
    superscriptaddress
]{revtex4-2}
\usepackage{hyperref}
\hypersetup{
    colorlinks=true,
    linkcolor=blue,
    filecolor=magenta,      
    urlcolor=cyan,
    pdftitle={Overleaf Example},
    pdfpagemode=FullScreen,
}
\usepackage{placeins}
\usepackage{amssymb}
\usepackage{float}
\usepackage{microtype}
\usepackage{amsmath}
\usepackage{amsthm}
\usepackage[english]{babel}
\usepackage{amsthm}
\usepackage{csquotes}
\theoremstyle{definition}
\newtheorem{definition}{Definition}[section]
\newtheorem{theorem}{Theorem}[section]
\newtheorem{lemma}[theorem]{Lemma}
\usepackage{geometry}
\usepackage{blindtext}
\usepackage{graphicx}
\usepackage{xcolor}   

\begin{document}

\title{Exotic newforms constructed from a linear combination of eta quotients}
\author{Anmol Kumar}
\email{anmolkumar@iisc.ac.in}
\affiliation{Indian Institute of Science, Bengaluru, India}

\date{\today}

\begin{abstract}
\sloppy
\raggedright
    K{\"o}hler, in [1], presented a weight 1 newform on $\Gamma_0(576)$ constructed from a linear combination of weight 1 eta quotients and asked, ``What would be a suitable $L$ and representation $\rho$ such that Deligne\text{-}Serre correspondence holds?" In this paper, we find the Galois field extension $L$ and representation $\rho$ such that the Deligne\text{-}Serre correspondence holds for this newform, and also study the splitting of primes in $L$ using the coefficients $a(p)$ of the newform. We also discuss an exotic newform on $\Gamma_0(1080)$ constructed from a linear combination of weight 1 eta quotients, find the corresponding Galois extension and representation, and study the splitting of primes in this extension. Furthermore, we find all such newforms that can be constructed from a linear combination of weight 1 eta quotients listed in [2] with $q$-expansion of the form $q+\sum_{k=2}^{\infty}a(k)q^k$.\\ %
\end{abstract}

\maketitle

\section{Introduction}
\begin{definition}
Let $\chi$ be a $\textit{dirichlet character modulo N}$ with $\chi(-1)=(-1)^k$. Let $$\Gamma_0(N):=\left\{\begin{pmatrix}a & b\\ c & d\end{pmatrix}\in \text{SL}_2(\mathbb{Z}):N\mid c   \right\}$$ A \textit{modular form of} \textit{type} $(k,\chi)$ on $\Gamma_0(N)$ is a function $f:\mathbb{H}\rightarrow\mathbb{C}$ such that 
\begin{itemize}
    \item $f$ is holomorphic on $\mathbb{H}$
    \item $f$ is holomorphic at all cusps of $\Gamma_0(N)$
    \item $f\left(\frac{az+b}{cz+d}\right)=\chi(d)(cz+d)^kf(z)$
\end{itemize}
    The vector space formed by modular forms of type $(k,\chi)$ on $\Gamma_0(N)$ is denoted by $M_k(\Gamma_0(N),\chi)$.
\end{definition}
\begin{definition}[Dedekind eta function]
The \textit{Dedekind eta function} is defined by $$\eta(z)=q^\frac{1}{24}\prod_{n=1}^{\infty} (1-q^n)$$
where $q=e^{2i\pi z}$, and $z$ lies in the upper half complex plane $\mathbb{H}$. 
\end{definition}
\begin{definition}
Define an \textit{eta-quotient} on level $N$ by the product $$ \prod_{0<m\mid N} \eta(mz)^{a_m}$$ where $a_m \in \mathbb{Z}$ $\forall$ $m\mid N$.
\end{definition}
\raggedright
The function $\eta(z)$ is a modular form of weight $\frac{1}{2}$ with a multiplier system on SL$_2(\mathbb{Z})$. In general, an eta-quotient $f(z)=\prod_{0<m\mid N} \eta(mz)^{a_m}$ is a meromorphic modular form of weight $k=\frac{1}{2}\sum_{m\mid N} a_m$ with a multiplier system on the congruence subgroup $\Gamma_0(N)$ for some $N$.
Yves Martin, in \cite{YM}, found all holomorphic eta-quotients $f$ of integral weight such that both $f(z)$ and its image under $\textit{Frickie involution}$ are eigenforms of all \textit{Hecke operators.} Bhattacharya, in \cite{Bha2}, gave a list of simple-holomorphic irreducible eta-quotients of weight 1 and also gave a simplified proof of Zagier's conjecture/Mersmann's theorem, which states that of any particular weight, there are only finitely many holomorphic eta quotients, none of which is an integral rescaling of another eta quotient or a product of two holomorphic eta quotients other than 1 and itself. Kohler in \cite{GK} constructed a linear combination of weight 1 eta quotients that was an exotic newform and asked the question, ``What would be a suitable $L$ and representation $\rho$ such that Deligne–Serre correspondence holds?'' In this paper, we determine $L$ and $\rho$ for this newform and study the splitting of primes in $L$ using the Fourier coefficients. We use the eta quotients of weight 1 listed in \cite{Bha1} to find a list of linear combinations of eta quotients of weight 1 that are newforms with q-expansion of the form $q+\sum_{k=2}^{\infty}a(k)q^k$, and also find out the projective image of the associated Galois representation. \\

\section{Prelimnaries}
We will need the following definitions and theorems.
\begin{theorem}\cite{Ono}
    If $f(z)=\prod_{0<m\mid N} \eta(mz)^{a_m}$ is an eta-quotient with $k=\frac{1}{2}\sum_{m\mid N} a_m \in \mathbb{Z},$ with the additional properties that $$\sum_{m\mid N}ma_m\equiv 0 \pmod {24}$$ and $$\sum_{m\mid N}\frac{N}{m}a_m\equiv 0 \pmod {24}$$, then $f(z)$ satisfies $$f\left(\frac{az+b}{cz+d}\right)=\chi(d)(cz+d)^kf(z)$$ for every $\begin{pmatrix}a & b\\ c & d\end{pmatrix}\in \Gamma_0(N).$
Here the character $\chi$ is a Dirichlet character modulo $N$ given by $$\chi(d)=\left(\frac{(-1)^ks}{d}\right)$$ where $s:=\displaystyle \prod_{d|N}d^{a_d}$. Moreover, if $f(z)$ is holomorphic at all cusps of $\Gamma_0(N)$, then $f(z)\in M_k(\Gamma_0(N),\chi)$.
\end{theorem}
\begin{definition}
    Let $\mathcal{L/K}$ be a Galois extension of number
    fields and denote the ring of integers of $\mathcal{K}$ and $\mathcal{L}$ by $\mathcal{O_K}$ and $\mathcal{O_L}$ respectively. Let $\mathfrak{p}$ be a prime in $\mathcal{O_K}$ and $\mathfrak{q}$ be a prime in $\mathcal{O_L}$ lying above $\mathfrak{p}$. Let $G=\operatorname{Gal}(\mathcal{L/K})$. Then, the $\textit{decomposition group}$ $D_{\mathfrak{q}}$ is defined as the stabilizer of $\mathfrak{q}$ in $G$.
\end{definition} The ramification index $e_\mathfrak{q}$ and residue field degree $f_\mathfrak{q}$ of $\mathfrak{q}$ are given by $e_\mathfrak{q}=\nu_\mathfrak{q}(\mathfrak{p}\mathcal{O_L})$ and $[\mathcal{O_L}/\mathfrak{q}:\mathcal{O_K/\mathfrak{p}}]$. Recall that $e_\mathfrak{q}$ and $f_\mathfrak{q}$ are independent of the choice of $\mathfrak{q}|\mathfrak{p}$. Hence, we can write $e_\mathfrak{p}$ and $f_\mathfrak{p}$ for $e_\mathfrak{q}$ and $f_\mathfrak{q}$ respectively for the rest of the paper.\\
\begin{lemma}\cite{Mar}
    The decomposition groups $D_\mathfrak{q}$ are all conjugates in $G$ for $\mathfrak{q}|\mathfrak{p}$ and $e_\mathfrak{p}f_\mathfrak{p}=|D_\mathfrak{q}|$. If $\mathfrak{p}$ is unramified in $\mathcal{K},$ then $e_\mathfrak{p}=1$, hence $f_p=|D_\mathfrak{q}|.$\\
\end{lemma}
Now, we want to determine the splitting behavior of an unramified prime $\mathfrak{p}$ of $\mathcal{O_K}$ in $\mathcal{O_L}$. Note that, by the $\textit{primitive element theorem,}$ there exists $\alpha$ in $\mathcal{O_L}$ such that $\mathcal{L}$ is generated over $\mathcal{K}$ by $\alpha$. Let the minimal polynomial of $\alpha$ over $\mathcal{O_K}$ be $h(x)$. We know that $\mathfrak{p}$ factorises in $\mathcal{O_L}$ as 
$$\mathfrak{p}\mathcal{O_L}=\prod_{i=1}^r\mathfrak{q}_i$$ where $r=|G|/f_\mathfrak{p}$ and $\mathfrak{q}_i$ are distinct prime ideals in $\mathcal{O_L}$ lying above $\mathfrak{p}$. Hence, we have $$h(x)\equiv\prod_{i=1}^rh_i(x)\text{ mod }\mathfrak{p}$$ where $h_i's$ are irreducible polynomials of degree $f_\mathfrak{p}$. Since $\mathfrak{p}$ is unramified, the map
$$D_\mathfrak{q}\rightarrow \operatorname{Gal}\big(\mathcal{O_L/\mathfrak{q}}\big{/}{\mathcal{O_K/\mathfrak{p}}}\big)$$
is an isomorphism of groups. The group $\operatorname{Gal}\big(\mathcal{O_L/\mathfrak{q}}\big{/}{\mathcal{O_K/\mathfrak{p}}}\big)$ is cyclic of order $f_\mathfrak{p}$. Therefore, $D_\mathfrak{q}$ is a cyclic group of order $f_\mathfrak{p}$ and is generated by the $\textit{Frobenius automorphism}$ $\phi_\mathfrak{q}$ defined as $$\phi_\mathfrak{q}:O_L/\mathfrak{q}\rightarrow O_L/\mathfrak{q}$$
$$x\mapsto x^{\mathcal{N(\mathfrak{p})}}$$ where $\mathcal{N(\mathfrak{p})}=|\mathcal{O_K}/\mathfrak{p}|$ is the norm of $\mathfrak{p}$. We define the $\textit{Frobenius element}$ $\sigma_\mathfrak{q}$ for an unramified prime $\mathfrak{q}|\mathfrak{p}$ to be the unique element in $G$ such that for all $x$ in $\mathcal{O_L}$, $\sigma_{\mathfrak{q}}(x)=x^{\mathcal{N}(\mathfrak{p})} \pmod{\mathfrak{q}}$. Note that if $\mathfrak{q}|\mathfrak{p}$ is unramified, then for all $\mathfrak{q'}|\mathfrak{p}$, $\sigma_{\mathfrak{q}}$ and $\sigma_{\mathfrak{q'}}$ are conjugates in $G$. Therefore, for $\mathfrak{q}|\mathfrak{p}$ unramified, we define the $\textit{Frobenius Class}$ of $\mathfrak{p}$ to be the conjugacy class of $\sigma_\mathfrak{q}$ in $G$ and denote any representative of the class by $\textrm{Frob}_\mathfrak{p}$.
\begin{definition}
Let $\mathcal{L/}\mathbb{Q}$ be a Galois extension with Galois group $G=\operatorname{Gal}(\mathcal{L/}\mathbb{Q})$ and a two-dimensional irreducible representation $\rho$. The $\textit{Artin L-function of $\rho$}$, denoted by $L(s,\rho)$, is defined as $$L(s,\rho)=\prod_pL_p(s,\rho)$$ where the local factors are given by 
\begin{align*}
    L_p(s,\rho)= \textrm{det}(I-\rho(\textrm{Frob}_p)p^{-s})^{-1} &=((1-\lambda_{1,p}p^{-s})(1-\lambda_{2,p}p^{-s}))^{-1} \\
    &=(1-\textrm{Tr}(\rho(\textrm{Frob}_p))p^{-s}+ \textrm{det}(\rho(\textrm{Frob}_p))p^{-2s})^{-1},
\end{align*}
where $\lambda_{1,p}$ and $\lambda_{2,p}$ are eigenvalues of $\textrm{Frob}_p$. Note that for ramified primes $p$, $\textrm{det}(\rho(\textrm{Frob}_p))=0$. \\
\end{definition}
\begin{theorem}[Chebotarev Density Theorem]\cite{Lag}
Let \( C \) be a conjugacy class of the Galois group \( G \) of a Galois extension \( K/\mathbb{Q} \). Then the density of the prime ideals \(\mathfrak{p} \in \mathcal{P}\) that are unramified in \( K \) and satisfy \(\mathrm{Frob}_{\mathfrak{p}} \in C\) is given by
\[
\delta\left(\{\mathfrak{p} \in \mathcal{P} \mid \mathfrak{p} \text{ is unramified and } \mathrm{Frob}_{\mathfrak{p}} \in C\}\right) = \frac{|C|}{|G|}.
\]
\end{theorem}
\begin{theorem}[Deligne-Serre]
    Let $f(z)=\sum_{n=1}^{\infty}a(n)q^n$ with $a(1)=1$ be a Hecke newform of weight 1 on $\Gamma_0(N)$ with character $\chi$ modulo $N$ satisfying $\chi(-1)=-1$. Assume that $f$ is an eigenform for all $T_p$ with $p\nmid N$, with eigenvalue $a_p$. Then there exists a Galois extension $L$ of $\mathbb{Q}$ and a 2-dimensional irreducible representation $\rho:\operatorname{Gal}{(L/\mathbb{Q})}\rightarrow GL_2(\mathbb{C})$ such that $\textrm{Tr}(\rho_f(\textrm{Frob}_p))=a_p$ and $\textrm{det}(\rho_f(\textrm{Frob}_p))=\chi(p)$ for all $p\nmid N$. Equivalently, there exists a Galois extension $L$ of $\mathbb{Q}$ and a 2-dimensional irreducible representation $\rho$ of G such that $L(s,\rho)=L_f(s)=\sum_{n=1}a(n)q^n$.
\end{theorem}
\begin{theorem}[Weil-Langlands-Khare-Wintenberger]
Given $\rho: G_\mathbb{Q}\rightarrow GL_2(\mathbb{C})$ be an irreducible, continuous, odd representation with Artin conductor $N$ and determinant $\chi$, let $L(s,\rho)=\sum_{n=1}\frac{a(n)}{n^s}$ be its Artin $L$-function. Then $f=\sum_{n=0}a(n)q^n$ is a normalized newform lying in $S_1(N,\chi)$.
\end{theorem}
\begin{definition}[Sturm Bound] For any space $M_k(\Gamma_0{(N)},\chi)$ of modular forms of weight $k$, level $N$, and character $\chi$, the Sturm bound is the integer
$$
B(M_k(\Gamma_0(N),\chi)) := \left\lfloor \frac{km}{12}\right\rfloor,$$
where
$$
m:=[SL_2(\mathbb{Z}):\Gamma_0(N)]=N\prod_{p|N}\left(1+\frac{1}{p}\right).
$$
\end{definition}
\begin{theorem}\cite{LMFDB} If $f=\sum_{n\ge 0}a_n q^n$ and $g=\sum_{n\ge 0}b_n q^n$ are elements of $M_k(N,\chi)$  with $a_n=b_n$ for all $n\le B(M_k(N,\chi))$ then $f=g$.
\end{theorem}
\begin{theorem}\cite{Lag} Let $f(z)=\sum_{n=1}^{\infty}a(n)q^n$ with $a(1)=1$ be a Hecke newform of weight 1 on $\Gamma_0(N)$ with character $\chi$ modulo $N$ satisfying $\chi(-1)=-1.$ Assume that $f$ is an eigenform for all $T_p$ with $p\nmid N$, with eigenvalue $a_p$. Let $L$ be the Galois number field such that the Deligne-Serre correspondence holds. Let $c(\rho(\textrm{Frob}_p)):=a_p^2/\chi(p)$ Then, the projective image of $G=\operatorname{Gal}(L/\mathbb{Q})$ is 
 \begin{itemize}
    \item $A_4$ if the proportion of primes $p$ with  $c(\rho(\textrm{Frob}_p))$ equal to 1, 2, and 3 converges to $\frac{1}{12}$, $\frac{1}{4}$, and $\frac{2}{3}$ respectively.
    \item $S_4$ if the proportion of primes $p$ with  $c(\rho(\textrm{Frob}_p))$ equal to 1, 2, 3 and 4 converges to $\frac{1}{25}$, $\frac{3}{8}$, $\frac{1}{3}$, and $\frac{1}{4}$ respectively.
    \item $A_5$ if the proportion of primes $p$ with  $c(\rho(\textrm{Frob}_p))$ equal to 1, 2, 3 and 5 converges to $\frac{1}{60}$, $\frac{1}{4}$, $\frac{1}{3}$, and $\frac{2}{5}$ respectively.
\end{itemize}
\end{theorem}
\section{Constructing newforms and the corresponding Number Fields}
Let $f(z)=\sum_{n=1}^\infty a(n)q^n$ with $a(1)=1$ be an eta quotient in $M_1(\Gamma_0(N),\chi)$, where $\chi(-1)=-1$. We compute the action of Hecke operator $T_{p_l}$ on $f$ for primes $p_l$ less than and coprime to $N$ and take a complex linear combination $\sum_{l}c_lT_{p_l}(f)$ of all linearly independent $T_{p_l}(f)$ which are eta quotients. The coefficients $c_l$ are chosen such that $\sum_{l}c_lT_{p_l}(f)$ is an eigenform for all $T_{p_l}$, where $p_l$ ranges through primes less than and coprime to N. This gives us a weight 1 modular form $F(z)=\sum_{n=1}^\infty b(n)q^n$ with character $\chi$ on $\Gamma_0(N)$. We calculate the approximate proportion of primes $p$ such that $b(p)=2$ by considering primes up to 10000. By $\textbf{Theorem 2.3}$, the reciprocal of this proportion gives us the order of group $G$. Next, we search for 2-dimensional complex irreducible representations of groups with this order, for which the fraction of elements with trace equal to $x$ is the same as the limiting proportion of primes $p$ with $b(p)=x$. Now, the primes $p$ such that $b(p)=2$ split completely in $L$; hence they split completely in all subfields of $L$. Using the \cite{LMFDB} website, we determine all subfields of $L$ and, hence, the composite field $L$. Using the fundamental theorem of Galois theory, we determine the subgroup structure of $G$, and hence the group $G$. We calculate the Frobenius element of prime $p$ up to the Sturm bound using \textit{sagemath}. Now, we find the Artin-$L$ function $L(s,\rho)$ for the representation $\rho$ and the $L$-function $L_F(s)=\sum_{n=1}\frac{b(n)}{n^s}$ up to the Sturm bound $B(M_k(\Gamma_0(N),\chi))$. $F$ is a newform if $L_F(s)=L(s,\rho)$, in which case the number field extension and representation such that the Deligne-Serre correspondence holds are $L$ and $\rho$, respectively. 
\section{Level 576}
In this section, we discuss a newform on $\Gamma_0(576)$ given in \cite{GK}. 
Consider the following $\eta$ quotients:
\begin{align*}
    f_{1}^{576}(z) &= \frac{\eta(4z)\eta(6z)^2\eta(12z)}{\eta(3z)\eta(24z)}, &
    f_{5}^{576}(z) &= \frac{\eta(2z)\eta(6z)\eta(12z)^2}{\eta(3z)\eta(24z)}, \\
    f_{7}^{576}(z) &= \frac{\eta(2z)\eta(4z)^2\eta(6z)}{\eta(z)\eta(8z)}, &
    f_{11}^{576}(z) &= \frac{\eta(2z)^2\eta(4z)\eta(12z)}{\eta(z)\eta(8z)}, \\
    f_{13}^{576}(z) &= \frac{\eta(z)\eta(8z)\eta(6z)}{\eta(2z)}, &
    f_{17}^{576}(z) &= \frac{\eta(z)\eta(8z)\eta(12z)}{\eta(4z)}, \\
    f_{19}^{576}(z) &= \frac{\eta(4z)\eta(3z)\eta(24z)}{\eta(12z)}, &
    f_{23}^{576}(z) &= \frac{\eta(2z)\eta(3z)\eta(24z)}{\eta(6z)}.
\end{align*}
The function $$F^{576}(24z)=f_{1}^{576}(24z)+f_{5}^{576}(24z)+f_{7}^{576}(24z)+f_{11}^{576}(24z)+i\sqrt{2}f_{13}^{576}(24z)-i\sqrt{2}f_{17}^{576}(24z)-i\sqrt{2}f_{19}^{576}(24z)+i\sqrt{2}f_{23}^{576}(24z)$$ defines a modular form of weight 1. The character is given by $\chi(d)=\Big(\frac{-24}{d}\Big)$. Consider the field extensions $K_1=\mathbb{Q}(\alpha)$, $K_2=\mathbb{Q}(\beta)$ and $K_3=\mathbb{Q}(\delta)$ of $\mathbb{Q}$ given by the defining polynomials $x^2 + 6$, $x^3 + 3x - 2$,  and $x^8 + 12x^6 + 42x^4 + 88x^2 - 6$ respectively. Let $L$ denote the composite field of the number fields $K_1$, $K_2$ and $K_3$. The Galois group $G=\operatorname{Gal}(L/\mathbb{Q})$ is $GL_2(\mathbb{F}_3)$. Now, consider the two-dimensional irreducible representation $\rho$ of $G$ 
$$\rho:G\rightarrow SL_2(\mathbb{C})$$ 
\FloatBarrier
\[
\begin{aligned}
&\rho\left(\begin{pmatrix} 1 & 0 \\ 0 & -1 \end{pmatrix}\right) =
\begin{pmatrix} 1 & 0 \\ 0 & -1 \end{pmatrix}, 
&&\rho\left(\begin{pmatrix} 1 & 1 \\ 0 & 1 \end{pmatrix}\right) =
\begin{pmatrix} \omega & 0 \\ 0 & \omega^2 \end{pmatrix}, \\
&\rho\left(\begin{pmatrix} 0 & 1 \\ -1 & 0 \end{pmatrix}\right) =
\begin{pmatrix} i & 0 \\ 0 & -i \end{pmatrix}, 
&&\rho\left(\begin{pmatrix} -1 & 1 \\ 0 & -1 \end{pmatrix}\right) =
\begin{pmatrix} -\omega & 0 \\ 0 & -\omega^2 \end{pmatrix}, \\
&\rho\left(\begin{pmatrix} 0 & 1 \\ 1 & -1 \end{pmatrix}\right) =
\begin{pmatrix} e^{\frac{i\pi}{4}} & 0 \\ 0 & -e^{\frac{-i\pi}{4}} \end{pmatrix}, 
&&\rho\left(\begin{pmatrix} 0 & 1 \\ 1 & 1 \end{pmatrix}\right) =
\begin{pmatrix} e^{\frac{-i\pi}{4}} & 0 \\ 0 & -e^{\frac{i\pi}{4}} \end{pmatrix}.
\end{aligned}
\]
\begin{table}[h]
$$
\begin{array}{c|rrrrrrrr}
  \rm \text{Class rep.}&\rm\big(\begin{smallmatrix} 1 & 0\\ 0 & 1 \end{smallmatrix}\big)&\rm\big(\begin{smallmatrix} -1 & 0\\ 0 & -1 \end{smallmatrix}\big)&\rm\big(\begin{smallmatrix} 1 & 0\\ 0 & -1 \end{smallmatrix}\big)&\rm\big(\begin{smallmatrix} 1 & 1\\ 0 & 1 \end{smallmatrix}\big)&\rm\big(\begin{smallmatrix} 0 & 1\\ -1 & 0 \end{smallmatrix}\big)&\rm\big(\begin{smallmatrix} -1 & 1\\ 0 & -1 \end{smallmatrix}\big)&\rm\big(\begin{smallmatrix}0 & 1\\ 1 & -1 \end{smallmatrix}\big)&\rm\big(\begin{smallmatrix}0 & 1\\ 1 & 1 \end{smallmatrix}\big)\cr
  \rm size&1&1&12&8&6&8&6&6\cr
\hline
  Tr(\rho)&2&-2&0&-1&0&1&-i\sqrt{2}&i\sqrt{2}\cr  
\end{array}
$$
\caption{Characters of the representation $\rho$}
\end{table}
\FloatBarrier
By $\textbf{Theorem 2.4}$, we know that there exists a Hecke eigenform $g(z)=\sum_{k=0}r(k)q^k$ with character $\chi'$ such that $\textrm{Tr}(\rho(\textrm{Frob}_p))=r(p)$ and $\textrm{det}(\rho(\textrm{Frob}_p))=\chi'(p)$. 
For primes $p<100$, we determined $\textrm{Frob}_p$ using \textit{sagemath} and the following behaviour was observed: 
\begin{table}[H]
\centering
 \begin{tabular}{||c c c c||} 
 \hline
 $a_p$ & Frobenius class representative of $p$ & $f_p$ & $L_p(s,\rho)^{-1}$ \\ [0.5ex] 
 \hline\hline
 2 & $\big(\begin{smallmatrix} 1 & 0\\ 0 & 1 \end{smallmatrix}\big)$ & 1 & $1-2p^{-s}+p^{-2s}$ \\ 
 -2 & $\big(\begin{smallmatrix} -1 & 0\\ 0 & -1 \end{smallmatrix}\big)$ & 2 & $1+2p^{-s}+p^{-2s}$ \\
 1 & $\big(\begin{smallmatrix} -1 & 1\\ 0 & -1 \end{smallmatrix}\big)$ & 6 & $1-p^{-s}+p^{-2s}$  \\
 -1 & $\big(\begin{smallmatrix} 1 & 1\\ 0 & 1 \end{smallmatrix}\big)$ & 3 & $1+p^{-s}+p^{-2s}$  \\
 0 & $\big(\begin{smallmatrix} 1 & 0\\ 0 & -1 \end{smallmatrix}\big)$ & 2 & $1-p^{-2s}$  \\
 0 & $\big(\begin{smallmatrix} 0 & 1\\ -1 & 0 \end{smallmatrix}\big)$ & 4 & $1+p^{-2s}$  \\
 $i\sqrt{2}$ & $\big(\begin{smallmatrix} 0 & 1\\ 1 & -1 \end{smallmatrix}\big)$ & 8 & $1-i\sqrt{2}p^{-s}-p^{-2s}$  \\
 $-i\sqrt{2}$ & $\big(\begin{smallmatrix} 0 & 1\\ 1 & 1 \end{smallmatrix}\big)$ & 8 & $1+i\sqrt{2}p^{-s}-p^{-2s}$  \\[1ex]
 \hline
 \end{tabular}
 \caption{}
\end{table}
Computing the first 100 terms of $L(s,\rho),$ we observe that $a(n)=r(n) \text{ }\forall \text{ }n<100$. Since $B(M_1(\Gamma_0(576),\chi))=96,$ by $\textbf{Theorem 2.5}$ we have $F^{576}=g$. Hence, $F^{576}$ is a newform and $L$ is the unique field extension of $\mathbb{Q}$ such that the Deligne-Serre correspondence holds. Using $\textbf{Theorem 2.6}$ and analyzing the data in $\textit{Python},$ we see that $F^{576}$ is an $S_4$ form.\\
\begin{theorem} [Splitting of Primes]
We present the splitting of primes in $L$ for primes $p$ not dividing the discriminant of $L$ in the following table: 
\FloatBarrier
\begin{table}[H]
\centering
 \begin{tabular}{||c c c c||} 
 \hline
 $a_p$ & Frobenius class representative of $p$ & $f_p$ & Number of primes $p$ splits into in $L$ \\ [0.5ex] 
 \hline\hline
 2 & $\big(\begin{smallmatrix} 1 & 0\\ 0 & 1 \end{smallmatrix}\big)$ & 1 & 48 \\ 
 -2 & $\big(\begin{smallmatrix} -1 & 0\\ 0 & -1 \end{smallmatrix}\big)$ & 2 & 24 \\
 1 & $\big(\begin{smallmatrix} -1 & 1\\ 0 & -1 \end{smallmatrix}\big)$ & 6 & 8  \\
 -1 & $\big(\begin{smallmatrix} 1 & 1\\ 0 & 1 \end{smallmatrix}\big)$ & 3 & 16  \\
 0 & $\big(\begin{smallmatrix} 1 & 0\\ 0 & -1 \end{smallmatrix}\big)$ & 2 & 24  \\
 0 & $\big(\begin{smallmatrix} 0 & 1\\ -1 & 0 \end{smallmatrix}\big)$ & 4 & 12  \\
 $i\sqrt{2}$ & $\big(\begin{smallmatrix} 0 & 1\\ 1 & -1 \end{smallmatrix}\big)$ & 8 & 6  \\
 $-i\sqrt{2}$ & $\big(\begin{smallmatrix} 0 & 1\\ 1 & 1 \end{smallmatrix}\big)$ & 8 & 6  \\[1ex]
 \hline
 \end{tabular}
 \caption{}
\end{table}
\FloatBarrier
\end{theorem}
\section{Level 1152}
In this section, we discuss a newform on $\Gamma_0(1152).$
$$f_{1}^{1152}(z)=\frac{\eta(144z)\eta(24z)^7\eta(16z)}{\eta(72z)^2\eta(8z)^2\eta(48z)^3}$$
$$f_{2}^{1152}(z)=\frac{T_{17}(f_{1}^{1152}(z))}{4}=\frac{\eta(72z)\eta(48z)^7\eta(8z)}{\eta(144z)^2\eta(24z)^3\eta(16z)^2}$$
The function
$$
F^{1152}(z)=\frac{f_{1}^{1152}(z)+2f_{2}^{1152}(z)}{3}
$$
defines a modular form of weight 1 on $\Gamma_0(1152)$. The character is given by
$$
\chi(d)=\Big(\frac{-2}{d}\Big).
$$
Consider the field extension \( L \) defined by the polynomial 
$$
p(x)=x^{8} - 4x^{7} +8x^{5} + 14x^{4} - 32x^{3} + 28x^{2} - 48x + 34.
$$
The Galois group \( G=\mathrm{Gal}(L/\mathbb{Q}) \) is the Quaternion group \( Q_8 = \langle a, b \mid a^4=1, b^2=a^2, bab^{-1}=a^{-1} \rangle \).
\FloatBarrier
\begin{table}[H]
\centering
\renewcommand{\arraystretch}{1.2} 
$$
\begin{array}{c|ccccc}
  \text{Class} & 1 & a^2 & a & b & ab \\ 
  \text{Size} & 1 & 1 & 2 & 2 & 2 \\ \hline
  \rho & 2 & -2 & 0 & 0 & 0 \\
\end{array}
$$
\caption{Characters of a 2-dimensional irreducible representation of $Q_8$}
\end{table}
\FloatBarrier
Using the method of analysis from \textit{section 4,} we see that $F^{1152}$ is a newform on $\Gamma_0(1152)$ and $L$ is the field extension such that the Deligne-Serre Correspondence holds. Moreover, $F^{1152}$ is a \textit{Dihedral form}.
\section{Level 5760}
$$
f_1^{5760}(z) = \frac{\eta(120z)^4 \eta(48z)}{\eta(240z)^2 \eta(24z)}
$$
$$
f_2^{5760}(z) = T_{29}(f_1^{240}(z)) = \frac{\eta(240z)^4 \eta(24z)}{\eta(48z) \eta(120z)^2}
$$
The function 
$$
F^{5760}(z) = f_1^{5760}(z) + 2i f_2^{5760}(z)
$$ 
defines a weight 1 modular form on \(\Gamma_0(5760)\). The character is 
$$
\chi(d) = \Big(\frac{-2}{d}\Big).
$$
Let \( L \) be the composite field of \( \mathbb{Q}[i] \) and the field extension of $\mathbb{Q}$ given by the defining polynomial 
$$
x^8 - 18x^4 + 9.
$$
The Galois group \( G = \operatorname{Gal}(L/\mathbb{Q}) \) is the Pauli group on 1-qubit:
$$
G = \langle a, b, c \mid a^4 = c^2 = 1, \, b^2 = a^2, \, ab = ba, \, ac = ca, \, cbc = a^2 b \rangle.
$$
\FloatBarrier
\begin{table}[H]
\centering
\renewcommand{\arraystretch}{1.3} 
$$
\begin{array}{c|cccccccccc}
  \text{Class} & 1 & a^2 & \text{$c$} & \text{$b^2$} & \text{$ab$} & a & a^3 & \text{$b$} & \text{$ac$} & \text{$abc$} \\ 
  \text{Size} & 1 & 1 & 2 & 2 & 2 & 1 & 1 & 2 & 2 & 2 \\ \hline
  \rho & 2 & -2 & 0 & 0 & 0 & -2i & 2i & 0 & 0 & 0 \\
\end{array}
$$
\caption{Characters of a 2-dimensional irreducible representation of \( G \)}
\end{table}
\FloatBarrier
Using the method of analysis from \textit{section 4,} we see that $F^{5760}$ is a newform on $\Gamma_0(5760)$ and $L$ is the field extension such that the Deligne-Serre Correspondence holds. Moreover, $F^{5760}$ is a $\textit{Dihedral form.}$
\section{Level 1080}
\[
f_1^{1080}(z) = \frac{\eta(36z)^2 \eta(60z) \eta(90z)}{\eta(18z) \eta(180z)}, \quad
f_2^{1080}(z) = \frac{\eta(12z)^2 \eta(30z) \eta(180z)}{\eta(60z) \eta(6z)}
\]
\[
f_3^{1080}(z) = \frac{\eta(12z) \eta(18z) \eta(180z)^2}{\eta(36z) \eta(90z)}, \quad
f_4^{1080}(z) = \frac{\eta(6z) \eta(36z) \eta(60z)^2}{\eta(12z) \eta(30z)}
\]
The function 
\[
F^{1080}(z) = f_1^{1080}(z) + i f_2^{1080}(z) + \zeta_8 f_3^{1080}(z) + i \zeta_8 f_4^{1080}(z),
\]
where \(\zeta_8 = \frac{\sqrt{2}i}{1-i}\), defines a weight 1 modular form on \(\Gamma_0(1080)\). The character is given by \(\chi(d) = \left(\frac{-15}{d}\right)\).
Consider the field extensions $K_1$, $K_2$, $K_3$, $K_4$ and $K_5$ of $\mathbb{Q}$ given by the defining polynomials $x^2-x+1$, $x^2-x-1$, $x^3-2$, $x^4-2x^3-4x-1$ and $x^{16}-4x^{14}-2x^{13}+10x^{12}-16x^{11}-4x^{10}+52x^9-41x^8+2x^7+80x^6-50x^5+8x^4+38x^3-20x^2+6x-1$ respectively. Let $L$ be the composite field of $K_1$, $K_2$, $K_3$, $K_4$ and $K_5$. 
The Galois group \( G = \operatorname{Gal}(L/\mathbb{Q}) \) is a central extension by \( C_4 \) of \( S_4 \):
\begin{multline*}
G = \langle a, b, c, d, e \mid a^4 = d^3 = e^2 = 1, \, b^2 = c^2 = a^2, \, ab = ba, \, ac = ca, \, ad = da, \, ae = ea, \\
cbc^{-1} = a^2b, \, dbd^{-1} = a^2bc, \, ebe = bc, \, dcd^{-1} = b, \, ece = a^2c, \, ede = d^{-1} \rangle.
\end{multline*}
\FloatBarrier
\begin{center}
\begin{small}
\begin{table}[!ht]
$
\begin{array}{c|rrrrrrrrrrrrrrrr}
  \rm \text{Class rep.}& 1 & a^2& ab & e&d& a& a^3& b& ae& a^2bcd& abce& bcd^2e& abcd^2e& bce& a^3bcd& ad\cr
  \rm size&1&1&6&12&8&1&1&6&12&8&6&6&6&6&8&8\cr
\hline
  \textrm{Tr}(\rho)&2&-2&0&0&-1&2i&-2i&0&0&1&-\sqrt{2}&-\sqrt{-2}&\sqrt{2}&\sqrt{-2}&-i&i\cr
\end{array}
$
\caption{Character table of $G$}
\end{table}
\end{small}
\end{center}
\FloatBarrier
Using the method of analysis from \textit{section 4,} we see that $F^{1080}$ is a newform on $\Gamma_0(1080)$ and $L$ is the field extension such that the Deligne-Serre Correspondence holds. Moreover, $F^{1080}$ is a $S_4$ \textit{form}.
\begin{theorem}
We present the splitting of primes in $L$ for primes $p$ not dividing the discriminant of $L$ in the following table: 
\FloatBarrier
\begin{table}[H]
\centering
 \begin{tabular}{||c c c c||} 
 \hline
 $a_p$ & Frobenius class representative of $p$ & $f_p$ & Number of primes $p$ splits into in $L$\\ [0.5ex] 
 \hline\hline
 2 & $1$ & 1 & $96$ \\ 
 -2 & $a^2$ & 2 & $48$ \\
 0 & $ab$ & 2 & $48$ \\ 
 0 & $e$ & 2 & $48$ \\ 
 0 & $b$ & 3 & $32$ \\ 
 0 & $ae$ & 4 & $24$ \\ 
 $2i$ & $a$ & 4 & $24$ \\ 
 $-2i$ & $a^3$ & 4 & $24$ \\ 
 1 & $a^2bcd$ & 4 & $24$ \\ 
 -1 & $d$ & 6 & $16$  \\
 $i$ & $ad$ & 8 & $12$  \\
 $-i$ & $a^3bcd$ & 8 & $12$  \\
 $\sqrt{2}$ & $abcd^2e$ & 8 & $12$ \\ 
 $-\sqrt{2}$ & $abce $ & 8 & $12$  \\
 $i\sqrt{2}$ & $bce$ & 12 & $8$  \\
 $-i\sqrt{2}$ & $bcd^2e$ & 12 & $8$  \\[1ex] 
 \hline
 \end{tabular}
 \caption{}
\end{table}
\FloatBarrier
\end{theorem}
\section{Level 9216}
$$f_1^{9216}(z)=\frac{\eta(192z)\eta(96z)^2\eta(48z)}{\eta(384z)\eta(24z)} \text{ , } f_2^{9216}(z)=\frac{T_5(f_1^{9216}(z))}{2}=\frac{\eta(48z)^3\eta(192z)^3}{\eta(24z)\eta(96z)^2\eta(384z)}$$
$$f_3^{9216}(z)=\frac{T_{13}(f_1^{9216}(z))}{2}=\frac{\eta(24z)\eta(96z)^4\eta(384z)}{\eta(192z)^2\eta(48z)^2} \text{ , } f_4^{9216}(z)=\frac{T_{17}(f_1^{9216}(z))}{4}=\eta(24z)\eta(384z)$$
The function $$F^{9216}(z)=f_1^{384}(z)+\sqrt{2}f_2^{384}(z)+\sqrt{2}f_3^{384}(z)-2f_4^{384}(z)$$
defines a modular form of weight 1 on $\Gamma_0(9216)$ and the character $\chi(d)=\big(\frac{-1}{d}$\big).
Consider the field extensions $K_1$ and $K_2$ of $\mathbb{Q}$ given by the defining polynomials $x^4+1$ and $x^8+72$, respectively. Let $L$ be the composite field of $K_1$ and $K_2$. 
The Galois group \( G = \operatorname{Gal}(L/\mathbb{Q}) \) is \( D_{16} \):
\[
G = \langle a, b \mid a^8 = b^2 = 1, \, bab = a^{-1} \rangle.
\]
\FloatBarrier
\begin{table}[H]
\centering
$$
\begin{array}{c|rrrrrrr}
  \text{Class} & \text{1} & \text{$a^4$} & \text{$b$} &\text{$ab$} & \text{$a^2$} & \text{$a$} & \text{$a^3$} \\
  \text{Size} & 1 & 1 & 4 & 4 & 2 & 2 & 2 \\
\hline
  \rho & 2 & -2 & 0 & 0 & 0 & -\sqrt{2} & \sqrt{2} \\
\end{array}
$$
\caption{Characters of a 2-dimensional irreducible representation of \( D_{16} \)}
\end{table}
\FloatBarrier
Using the method of analysis from \textit{section 4,} we see that $F^{9216}$ is a newform on $\Gamma_0(9216)$ and $L$ is the field extension such that the Deligne-Serre Correspondence holds. Moreover, $F^{9216}$ is a $\textit{Dihedral form.}$
\section{Level 23040}
\raggedright
$$f_1^{23040}(z)=\frac{\eta(120z)^2\eta(96z)^2\eta(48z)}{\eta(240z)\eta(192z)\eta(24z)}
\text{ , } f_2^{23040}(z)=\frac{T_7(f_1^{23040}(z))}{2}=\frac{\eta(240z)^2\eta(192z)\eta(48z)^2}{\eta(480z)\eta(96z)\eta(24z)}$$
$$f_3^{23040}(z)=\frac{T_{31}(f_1^{23040}(z))}{2}=\frac{\eta(24z)\eta(192z)\eta(240z)^5}{\eta(48z)\eta(120z)^2\eta(480z)^2}\text{ , } f_{4}^{23040}(z)=-\frac{T_{29}(f_1^{23040}(z))}{2}=\frac{\eta(24z)\eta(96z)^2\eta(480z)^2}{\eta(48z)\eta(192z)\eta(240z)}$$
$$f_{5}^{23040}(z)=-\frac{T_{53}(f_1^{23040}(z))}{2}=\frac{\eta(192z)\eta(24z)\eta(480z)^5}{\eta(96z)\eta(240)^2\eta(960z)^2}\text{ , } f_{6}^{23040}(z)=\frac{T_{59}(f_1^{23040}(z))}{4}=\frac{\eta(480z)^2\eta(192z)\eta(48z)^3}{\eta(240z)\eta(24z)\eta(96z)^2}$$
$$f_{7}^{23040}(z)=\frac{T_{73}(f_1^{23040}(z))}{2}=\frac{\eta(24z)\eta(96z)^3\eta(240z)^2}{\eta(48z)^2\eta(192z)\eta(480z)}\text{ , } f_{8}^{23040}(z)=-\frac{T_{83}(f_1^{23040}(z))}{4}=\frac{\eta(48z)^2\eta(96z)\eta(960z)^2}{\eta(24z)\eta(192z)\eta(480z)}$$
\begin{align*}
F^{23040}(z) &= \frac{
f_1^{23040}(z) + \sqrt{2} f_2^{23040}(z) + i\sqrt{2} f_3^{23040}(z) + i\sqrt{2} f_4^{23040}(z)
}{1+i} \\
&\quad + \frac{
\sqrt{2} f_5^{23040}(z) + 2 f_6^{23040}(z) + i f_7^{23040}(z) + 2i f_8^{23040}(z)
}{1+i}
\end{align*}

The function $F^{23040}$ defines a weight 1 modular form of weight 1 on $\Gamma_0(23040)$. The character is $\chi(d)=\Big(\frac{-10}{d}\Big)$.
Consider the field extensions $K_1$ and $K_2$ of $\mathbf{Q}$ given by the defining polynomials $x^2+2$ and $x^{16} - 60x^{14} + 1098x^{12} - 8280x^{10} + 28710x^8 - 43200{x^6} + 19116x^{4} + 6480x^{2} + 324$ respectively. Let $L$ be the composite field of $K_1$ and $K_2$.
The Galois group \( G = \mathrm{Gal}(L/\mathbb{Q}) \) is the central product of \( C_4 \) and \( D_8 \):
\[
G = C_4 \circ D_8 = \langle a, b, c \mid a^4 = c^2 = 1, \, b^4 = a^2, \, ab = ba, \, ac = ca, \, cbc = a^2b^3 \rangle.
\]
\FloatBarrier
\begin{table}[H]
\centering
$$
\begin{array}{c|rrrrrrrrrrrrrr}
  \text{Class} & \text{1} & \text{$a^2$} & \text{$c$} & \text{$bc$} & \text{$ab^2$} & \text{$a$} & \text{$a^3$} & \text{$b^2$} & \text{$ac$} & \text{$abc$} & \text{$b$} & \text{$b^3$} & \text{$ab$}& \text{$ab^3$} \\
  \text{Size} & 1 & 1 & 2 & 4 & 4 & 1 & 1 & 2 & 4 & 4 & 2 & 2 & 2 & 2 \\
\hline
  \rho & 2 & -2 & 0 & 0 & 0 & -2i & 2i & 0 & 0 & 0 & -\sqrt{-2} & \sqrt{-2} & -\sqrt{2} & \sqrt{2} \\
\end{array}
$$
\caption{Characters of a 2-dimensional irreducible representation of \( C_4 \circ D_8 \)}
\end{table}
\FloatBarrier
Using the method of analysis from \textit{section 4,} we see that $F^{23040}$ is a newform on $\Gamma_0(23040)$ and $L$ is the field extension such that the Deligne-Serre Correspondence holds. Moreover, $F^{23040}$ is a $\textit{Dihedral form.}$
\section{Main result}
The analysis done in previous sections proves the result:
\begin{theorem}Within the weight 1 eta quotients listed in \cite{Bha1} but $\textbf{not}$ in \cite{GK2}, the following are all possible linear combinations that give a Hecke eigenform with $q$-expansion starting at $q$
\setlength{\tabcolsep}{12pt} 
\renewcommand{\arraystretch}{1.5} 
\begin{table}[H]
\centering
\begin{tabular}{|| c | c ||} 
 \hline
 \textbf{Newform} & \textbf{Projective Image of Galois Representation} \\ [1ex] 
 \hline\hline
 $F^{1152}$ & Dihedral \\ \hline
 $F^{576}$ & $S_4$ \\ \hline
 $F^{5760}$ & Dihedral \\ \hline
 $F^{1080}$ & $S_4$ \\ \hline
 $F^{9216}$ & Dihedral \\ \hline
 $F^{23040}$ & Dihedral \\ 
 \hline
\end{tabular}
\caption{Projective images of Galois representations corresponding to the newforms.}
\label{tab:proj_image}
\end{table}
\end{theorem}
\section{Conclusion}
We have discussed the splitting of primes in the corresponding Galois number fields using the coefficients $a(p)$ of the newforms. We also proved that the representation is of type $S_4$ for both the newforms. Bhattacharya, in \cite{Bha1}, lists simple holomorphic eta quotients on $\Gamma_0(N)$ of weight-1 up to level 768. Mersmann, in \cite{Mer}, shows that there are only finitely many simple holomorphic eta quotients of weight-1. Using the list of weight-1 simple holomorphic eta quotients, we attempted to find all exotic newforms that can be constructed from a linear combination of holomorphic eta quotients. Since there are only finitely many simple holomorphic eta quotients of weight-$1$, the number of such newforms is finite. The two newforms discussed in this paper are the only exotic newforms with the $q$-series being $q+\sum_{n=2}a(n)q^n$ that can be constructed from a linear combination of holomorphic eta quotients contained in this list. One can attempt to find a complete list of all exotic newforms of weight-1 that can be constructed from a linear combination of simple holomorphic eta quotients. The newforms discussed in the paper further raise a question: If an exotic newform is constructed from a linear combination of simple holomorphic eta quotients, is the representation $\rho$ such that the Deligne-Serre correspondence holds, always of type $S_4$?
\section{Acknowledgements}
I would like to express my deepest gratitude to Prof. Yingkun Li for encouraging and guiding me throughout the course of this research. The research was carried out at Technische Universität Darmstadt in Germany, and I am grateful to DAAD-WISE for the funding. Lastly, I would like to extend my heartfelt thanks to my family and friends for their constant support and motivation.

\bibliographystyle{apsrev4-1} 
\bibliography{references}          

\end{document}